\documentclass[12pt,reqno]{amsart}

\usepackage{mathrsfs}
\usepackage{amsmath}
\usepackage{amsfonts}
\usepackage{amssymb}
\usepackage[all,cmtip]{xy}
\usepackage{hyperref}
\usepackage{enumerate}  
\usepackage{enumitem}
\usepackage{ulem}
\usepackage[all]{xy}
\usepackage[toc,page]{appendix}
\usepackage{tikz-cd}
\usepackage[all,cmtip]{xy}
\usepackage[toc,page]{appendix}
\usepackage{amsmath, amsthm, amsfonts, amssymb}
\usepackage{graphicx}
\usepackage{subcaption}
\usepackage{tikz-cd}
\usepackage[utf8]{inputenc}
\usepackage[english]{babel}
\usepackage{dsfont}
\usepackage{xcolor}
\usepackage{arydshln}
\usepackage{faktor}
\usepackage{mathtools}
\usepackage{enumitem}

\usepackage{dsfont}

\makeatletter
\newcommand*{\da@rightarrow}{\mathchar"0\hexnumber@\symAMSa 4B }
\newcommand*{\da@leftarrow}{\mathchar"0\hexnumber@\symAMSa 4C }
\newcommand*{\xdashrightarrow}[2][]{%
	\mathrel{%
		\mathpalette{\da@xarrow{#1}{#2}{}\da@rightarrow{\,}{}}{}%
	}%
}
\newcommand{\xdashleftarrow}[2][]{%
	\mathrel{%
		\mathpalette{\da@xarrow{#1}{#2}\da@leftarrow{}{}{\,}}{}%
	}%
}
\newcommand{\xdashdownarrow}[2][]{%
	\mathrel{%
		\mathpalette{\da@xarrow{#1}{#2}\da@downarrow{}{}{\,}}{}%
	}%
}
\newcommand*{\da@xarrow}[7]{%
	\sbox0{$\ifx#7\scriptstyle\scriptscriptstyle\else\scriptstyle\fi#5#1#6\m@th$}%
	\sbox2{$\ifx#7\scriptstyle\scriptscriptstyle\else\scriptstyle\fi#5#2#6\m@th$}%
	\sbox4{$#7\dabar@\m@th$}%
	\dimen@=\wd0 %
	\ifdim\wd2 >\dimen@
	\dimen@=\wd2 %   
	\fi
	\count@=2 %
	\def\da@bars{\dabar@\dabar@}%
	\@whiledim\count@\wd4<\dimen@\do{%
		\advance\count@\@ne
		\expandafter\def\expandafter\da@bars\expandafter{%
			\da@bars
			\dabar@ 
		}%
	}%  
	\mathrel{#3}%
	\mathrel{%   
		\mathop{\da@bars}\limits
		\ifx\\#1\\%
		\else
		_{\copy0}%
		\fi
		\ifx\\#2\\%
		\else
		^{\copy2}%
		\fi
	}%   
	\mathrel{#4}%
}
\makeatother

\usepackage{graphicx}

\newtheorem{theorem}{Theorem}[section]

\theoremstyle{definition}

\newtheorem{remark}[theorem]{Remark}

\date{}



 \def\dashmapsto{\mathrel{\mapstochar\dashrightarrow}}

\begin{document}

\title[A note on representations of  $SL_n$]{\bf  A note on representations of  $SL_n$}
\author{Y\lowercase{isha} Y\lowercase{ao} 
   \ \\
Y\lowercase{ale} U\lowercase{niversity}\\ \ }

\vspace{3cm} \maketitle

\begin{abstract}
In this note, we study irreducible unitary representations of special linear groups of lower ranks, in terms of the matrix models of Gelfand-Naimark and Gelfand-Graev. Review of existing literature is provided. We also add some new calculation based on existing theory. 
\end{abstract}

\tableofcontents

\section{General theory}\label{basict}

%In this section, we shall give an explicit construction of the irreducible unitary representations of $SL_n(\mathbb F)$ for $n=3,4$ and $\mathbb F=\mathbb R,\mathbb C$. Before that, we first recall the corresponding results for $SL_2(\mathbb C)$ and $SL_2(\mathbb R)$.   

Gelfand-Naimark \cite{GN1} proved that the non-trivial irreducible unitary representations of $SL_2(\mathbb C)$ are given by the  principal and complementary series as follows. 

Define an upper triangular matrix group
\begin{equation}
	K:=\left\{\left. k=\left(
	\begin{matrix}
		k_{11}&k_{12}\\
		0&k_{22}\\
	\end{matrix}\right)\right\vert_{} 
	k_{ij}\in\mathbb C,\, k_{11}\cdot k_{22}=1 \right\}.
\end{equation}
For $m_2\in\mathbb Z$ and $\rho_2\in\mathbb R$, define a complex character $\chi$ of $K$ by
\begin{equation}\label{cc}
	\chi(k):=|k_{22}|^{m_2+\sqrt{-1}\rho_2}\cdot k_{22}^{-m_2}, \,\,k\in K.
\end{equation}
Write $\chi=(m_2,\rho_2)$ to indicate the parameters. Denote by $\mathcal H_{\chi}$ the Hilbert space of square integrable (complex-valued) functions on the complex plane with respect to the Lebesgue measure $\frac{i}{2}dz\wedge d\bar{z}$.
Now, the principal series $T^{\chi}$ of $SL_2(\mathbb C)$ can be described as follows. For each $g=(a_{ij})\in SL_2(\mathbb C)$,  define a linear operator $T^{\chi}(g)$ such that it acts on 
$\mathcal H_{\chi}$ by
\begin{equation}\label{ps}
	\left(T^{\chi}(g)f\right)(z):=|a_{12}z+a_{22}|^{-2+m_2+\sqrt{-1}\rho_2}\left(a_{12}z+a_{22}\right)^{-m_2}f\left(\frac{a_{11}z+a_{21}}{a_{12}z+a_{22}}\right),\,f\in\mathcal H_{\chi}.
\end{equation}
It is well-known that $T^{\chi}$, $T^{\chi^{\prime}}$ are equivalently if and only if $(m^{\prime}_2,\rho^{\prime}_2)=\pm(m_2,\rho_2)$. 

The complementary series  $T^{\chi}$ can be similarly given, based on the complex character $\chi=(\sigma)$ defined by $\chi(k):=|k_{22}|^{-2\sigma}$, $k\in K$, $0<\sigma<1$. 
Let $\mathcal H_{\chi}$ be the Hilbert space of square integrable functions on the complex plane with respect to the norm 
\begin{equation}
	(f,f):=\left(\frac{i}{2}\right)^2\iint_{\mathbb C\times\mathbb C}\frac{f(z)\cdot\overline{f(w)}}{|z-w|^{2-2\sigma}}\,dz\wedge d\bar z\wedge dw\wedge d\bar w.%,\,\, f,h\in\mathcal H_{\chi}.
\end{equation}
The transformation law is 
\begin{equation}
	\left(T^{\chi}(g)f\right)(z):=|a_{12}z+a_{22}|^{-2-2\sigma}f\left(\frac{a_{11}z+a_{21}}{a_{12}z+a_{22}}\right),\, f\in\mathcal H_{\chi}.
\end{equation}
Then $T^{\chi}$, $0<\sigma<1$, are unitary,  irreducible, and inequivalent.

Bargmann \cite{B} proved that the non-trivial irreducible unitary representations of $SL_2(\mathbb R)$ consist of the following principal,  complementary, discrete, and limits of discrete series. 

Let 
\begin{equation}
	K:=\left\{\left.  k=\left(
	\begin{matrix}
		k_{11}&k_{12}\\
		0&k_{22}\\
	\end{matrix}\right)\right\vert_{} 
	k_{ij}\in\mathbb R,\, k_{11}\cdot k_{22}=1 \right\}.
\end{equation}
For $m_2\in\{0,1\}$ and $\rho_2\in\mathbb R$, define a character  $\chi=(m_2,\rho_2)$ by
\begin{equation}\label{indu}
	\chi(k):=|k_{22}|^{\sqrt{-1} \rho_2} \cdot\left(\frac{k_{22}}{|k_{22}|}\right)^{m_2},\,\,k\in K.
\end{equation}
Let $\mathcal H_{\chi}$ the Hilbert space of square integrable (complex-valued) functions on the real line $(-\infty,\infty)$ with respect to the Lebesgue measure $dx$. The principal series $T^{\chi}$ can thus be realized as linear operators acting on $\mathcal H_{\chi}$ by
\begin{equation}\label{rps}
	\left(T^{\chi}(g)f\right)(x):=|a_{12}x+a_{22}|^{-1+\sqrt{-1}\rho_2}\left(\frac{a_{12}x+a_{22}}{|a_{12}x+a_{22}|}\right)^{m_2}f\left(\frac{a_{11}x+a_{21}}{a_{12}x+a_{22}}\right),\,\,f\in\mathcal  H_{\chi},
\end{equation}
where $g=(a_{ij})\in SL(2,\mathbb R)$. Notice that $T^{\chi}$, $T^{\chi^{\prime}}$ are  equivalently if and only if $m^{\prime}_2=m_2$ and $\rho^{\prime}_2=\pm\rho_2$; $T^{\chi}$ is irreducible if and only if $(m_2,\rho_2)\neq (1,0)$.

Realize the complementary series $T^{\chi}$ as follows. For $0<\sigma<1$, define $\chi(k):=|k_{22}|^{-\sigma}$, $k\in K$.  Let $\mathcal H_{\chi}$ be the Hilbert space  of square integrable functions on the real line with respect to the norm
\begin{equation}
	(f,f):=\frac{1}{\Gamma(\sigma)}\iint_{\mathbb R\times\mathbb R}\frac{f(x)\cdot\overline{f(y)}}{|x-y|^{1-\sigma}}\,dx\,dy.%,\,\,f,h\in\mathcal  H_{\chi}.
\end{equation}
For $g=(a_{ij})\in SL(2,\mathbb R)$, define
\begin{equation}
	\left(T^{\chi}(g)\right)f(x):=|a_{12}x+a_{22}|^{-1-\sigma}f\left(\frac{a_{11}x+a_{21}}{a_{12}x+a_{22}}\right),\,\, f\in \mathcal H_{\chi}.
\end{equation}
Then $T^{\chi}$, $0<\sigma<1$, are unitary, irreducible, and inequivalent.

The discrete series representation $T^{+,s}$ (resp. $T^{-,s}$),  $s=1,2,\cdots$, is defined on the Hilbert space $H^+_{s}$ (resp. $H^-_{s}$) of square integrable analytic functions on the upper half-plane (resp. lower half-plane) with respect to the norm
\begin{equation}\label{dinn}
	(f,f):=\frac{i}{2}\int_{\Im z>0} |f(z)|^2\left|\Im z\right|^{s-1}dz\wedge d\bar z\,\,%,\,f,h\in H^+_{s} 
	\left({\rm resp.}\,\,\,\frac{i}{2}\int_{\Im z<0} |f(z)|^2\left|\Im z\right|^{s-1}dz\wedge d\bar z%,\,f,h\in H^-_{s}
	\right).
\end{equation}
The transformation law is 
\begin{equation}\label{d0}
	\left(T^{\pm,s}(g)\right)f(z):=\left(a_{12}z+a_{22}\right)^{-1-s}f\left(\frac{a_{11}z+a_{21}}{a_{12}z+a_{22}}\right),\,\,f\in H^{\pm}_{s}.
\end{equation}
%It is well-known that $T^{\pm,s}$, $s=1,2,\cdots$, are unitary, irreducible, and inequivalent.

When $\chi=(m_2,\rho_2)=(1,0)$, the unitary representation $T^{\chi}$ defined by  (\ref{rps}) is no longer irreducible, but rather a direct sum of two irreducible subrepresentations $T^{+,0}$ and $T^{-,0}$. More precisely, $\mathcal H_{(1,0)}=H^2(\mathbb H^{+})\oplus H^2(\mathbb H^{-})$, where
$H^2(\mathbb H^{+})$ (resp. $H^2(\mathbb H^{-})$) is the Hardy space of analytic functions on the upper half-plane (resp. lower half-plane) with the norm
\begin{equation}
	\|f\|^2_{H^2(\mathbb H^{+})}:=\sup_{y>0}\int_{-\infty}^{\infty}|f(x+iy)|^2 dx\,\,\,\left({\rm resp.}\,\,\|f\|^2_{H^2(\mathbb H^{-})}:=\sup_{y<0}\int_{-\infty}^{\infty}|f(x+iy)|^2 dx\right).
\end{equation}
(See \cite{Dur} for more analytical aspects of Hardy spaces.) Then $T^{\pm,0}(g)$ acts on  $H^2(\mathbb H^{\pm})$ by
\begin{equation}\label{d2}
	\left(T^{\pm,s}(g)\right)f(z):=\left(a_{12}z+a_{22}\right)^{-1}f\left(\frac{a_{11}z+a_{21}}{a_{12}z+a_{22}}\right),\,\,f\in H^2(\mathbb H^{\pm}).
\end{equation} 
\medskip

The unitary duals of general linear groups over  Archimedean fields  are determined by \cite{Vo, Ta1,Ta2}.
Following \cite{Ta2}, we next recall certain general facts for the irreducible unitary representations of $GL_n(\mathbb F)$, when $\mathbb F=\mathbb R,\mathbb C$.

Given representations $\sigma_i$ of $GL_{n_i}(\mathbb F)$, $1\leq i\leq k$, $\sigma:= \sigma_1\otimes\sigma_2\cdots\otimes\sigma_k$ is a representation of $M:=GL_{n_1}(\mathbb F)\times\cdots\times GL_{n_k}(\mathbb F)$. Denote by ${\rm Ind}(\sigma)$ the representation of $GL_{n_1+\cdots+n_k}(\mathbb F)$ parabolically induced
by $\sigma$ from  the unique standard parabolic subgroup of $GL_{n_1+\cdots+n_k}(\mathbb F)$ associated with the Levi factor $M$. 

Let $\widehat{\mathbb F^{\times}}$ be the set of unitary characters of $\mathbb F^{\times}$. Let $\det_n$  be the determinant homomorphism of $GL_n(\mathbb F)$.   For
$\chi\in\widehat{\mathbb F^{\times}}$, $\chi\circ\det_n$ is a character of  $GL_n(\mathbb F)$.  Let $|\cdot|_{\mathbb F}$ the modulus character of $\mathbb F$; then, $|\cdot|_{\mathbb R}$ is the standard absolute value on $\mathbb R$, and
$|\cdot|_{\mathbb C}$ is the square of the standard absolute value on $\mathbb C$. When $\mathbb F=\mathbb C$, let
\begin{equation}
	B:=\left\{\chi\circ{\rm det}_{n},\,{\rm Ind}\left(|{\rm det}_{n}|_{\mathbb C}^{\alpha}(\chi\circ{\rm det}_{n})\otimes|{\rm det}_{n}|_{\mathbb C}^{-\alpha}(\chi\circ{\rm det}_{n})\right)\left|\begin{matrix}
	\,\chi\in\widehat{\mathbb C^{\times}},n\in\mathbb Z^+,\\0<\alpha<\frac{1}{2}
	\end{matrix}\right.\right\}.
\end{equation} 
Then,
\begin{theorem}[\cite{Ta2}]\label{glnc} 
For $\sigma_1,\cdots,\sigma_k\in B$, the parabolically induced representation 
\begin{equation}\label{ind}
{\rm Ind}(\sigma_1\otimes\sigma_2\cdots\otimes\sigma_k)
\end{equation}
is an irreducible unitarizable representation of a general linear group over $\mathbb C$.
Each irreducible unitarizable representation $\pi$ of $GL_n(\mathbb C)$ is equivalent to a parabolically induced representation from (\ref{ind}).  Moreover, $\pi$ determines the sequence $\sigma_1,\cdots,\sigma_k$ uniquely up to a permutation.  %Further, if $p$ is a permutation of $\{1,\cdots,k\}$, then \begin{equation}{\rm Ind}(\sigma_1\otimes\sigma_2\cdots\otimes\sigma_k)={\rm Ind}(\sigma_{p(1)}\otimes\sigma_{p(2)}\cdots\otimes\sigma_{p(k)}).\end{equation}.
\end{theorem}

\begin{remark}
Hence, according to Mackey analysis,  the irreducible unitary representations of $SL_n(\mathbb C)$ are exhausted by Gelfand-Naimark's construction \cite{GN2} together with Stein's complementary series \cite{St}.  
\end{remark}

For an irreducible representation $\delta$ of $GL_2(\mathbb R)$ which is square integrable modulo the center, and a positive integer $n$, the parabolically induced representation
\begin{equation}
	{\rm Ind} \left(|\det|_{\mathbb F}^{\frac{n-1}{2}}\delta\otimes|\det|_{\mathbb F}^{\frac{n-1}{2}-1}\delta\otimes\cdots\otimes|\det|_{\mathbb F} ^{-\frac{n-1}{2}}\delta\right)
\end{equation}
has a unique irreducible quotient, which will be denoted by $u(\delta,n)$. Define 
\begin{equation}
	\begin{split}
		B&:=\left\{{u(\delta,n),\rm Ind}\left(|{\rm det}_n|_{\mathbb R}^{\alpha}u(\delta,n)\otimes|{\rm det}_n|_{\mathbb R}^{-\alpha}u(\delta,n)\right)\left|\begin{matrix}
		\delta\in\widehat {GL_2(\mathbb R)} {\rm\,\,square\, integrable}\\
		    {\rm modulo\, the\, center,\,} n\in\mathbb Z^{+},\\\,0<\alpha<\frac{1}{2}\\
		\end{matrix}\right.\right\}\\
	&\,\,\,\,\,\,\,\,\,\,\,\,\,\,\,\bigcup\left\{\chi\circ{\rm det}_n,\,{\rm Ind}\left(|{\rm det}_n|_{\mathbb R}^{\alpha}(\chi\circ{\rm det}_n)\otimes|{\rm det}_n|_{\mathbb R}^{-\alpha}(\chi\circ{\rm det}_n)\right)\left|\begin{matrix}\chi\in\widehat{\mathbb R^{\times}},\,n\in\mathbb Z^{+}\\
		0<\alpha<\frac{1}{2}
	\end{matrix}\right.\right\}.\\
	\end{split}
\end{equation}
Then,
\begin{theorem}[\cite{Ta2}]\label{glnr} 
For $\sigma_1,\cdots,\sigma_k\in B$, the parabolically induced representation 
\begin{equation}\label{rind}
{\rm Ind}(\sigma_1\otimes\sigma_2\cdots\otimes\sigma_k)
\end{equation}
is an irreducible unitarizable representation of a general linear group over $\mathbb R$. Each irreducible unitarizable representation $\pi$ of $GL_n(\mathbb R)$ is equivalent to a parabolically induced representation from (\ref{rind}). Moreover, $\pi$ determines the sequence $\sigma_1,\cdots,\sigma_k$uniquely up to a
permutation.
\end{theorem}

\begin{remark}\label{split}
An irreducible unitary representation of $SL_n(\mathbb R)$ can be extended to an irreducible representation of $\mathbb R\times SL_n(\mathbb R)$ and then induced up to $GL_n(\mathbb R)$. If we  restrict back the induced representation of $GL_n(\mathbb R)$ to $SL_n(\mathbb R)$, we shall get a sum of $k$ irreducible representations depending on the index of $\mathbb R\times SL_n(\mathbb R)$ in $GL_n(\mathbb R)$ ($k=1$ when $n$ is odd and $k=2$ when $n$ is even). 
\end{remark}

Based on the above theorems, in the remaining of this note, we shall describe explicitly the unitary representations of low-rank special linear groups over the complex and real fields (see \cite{GN2,GG1,Ts,Va,Sp,Duf} for reference).

\section{Irreducible unitary representations of \texorpdfstring{$SL_3(\mathbb C)$} {rr}}\label{sl3c} Following \cite{GN2}, we denote by $K$ the upper triangular matrix subgroup of $SL_3(\mathbb C)$, which consists of complex matrices $k=(k_{pq})$ such that $k_{pq}=0$ for $p>q$. Denote by $Z$ the unipotent lower triangular matrix subgroup of $SL_3(\mathbb C)$, which consists of matrices $z=(z_{pq})$ such that $z_{pp}=1$ and $z_{pq}=0$ for $p<q$. Equip $K$ with a left Haar measure
\begin{equation}
	d\mu_l(k):=\left(\frac{i}{2}\right)^{5}|k_{33}|^2\cdot dk_{12}\wedge d\overline{k_{12}}\wedge dk_{13}\wedge d\overline{k_{13}}\wedge dk_{22}\wedge d\overline{k_{22}}\wedge dk_{23}\wedge d\overline{k_{23}}\wedge dk_{33}\wedge d\overline{k_{33}},
\end{equation}
and a right Haar measure 
\begin{equation}
	d\mu_r(k):=\left(\frac{i}{2}\right)^{5}|k_{22}|^{-4}\cdot|k_{33}|^{-6}\cdot dk_{12}\wedge d\overline{k_{12}}\wedge dk_{13}\wedge d\overline{k_{13}}\wedge dk_{22}\wedge d\overline{k_{22}}\wedge dk_{23}\wedge d\overline{k_{23}}\wedge dk_{33}\wedge d\overline{k_{33}}.
\end{equation}
Then the modular function $\beta$ of  $K$ is given by 
\begin{equation}\label{beta}
\beta(k)=\frac{d\mu_l(k)}{d\mu_r(k)}=|k_{22}|^{4}\cdot |k_{33}|^{8}.
\end{equation}

Since $Z$ is an open subset of the flag variety $SL_3(\mathbb C)/K$, each $g\in SL_3(\mathbb C)$ induces a birational transformation $\bar g:z\dashmapsto z\bar g$ of $Z$. More explicitly, for $g\in SL_3(\mathbb C)$ and generic $z\in Z$,
\begin{equation}\label{zzp}
	\left(\begin{matrix}
		1&0 &0\\
		z_{21}&1&0\\
		z_{31}&z_{32}&1\\
	\end{matrix}\,\right) \left(\begin{matrix}
		g_{11}&g_{12}&g_{13}\\
		g_{21}&g_{22}&g_{23}\\
		g_{31}&g_{32} &g_{33}\\
	\end{matrix}
	\right)=\left(\begin{matrix}
		k_{11}&k_{12} &k_{13}\\
		0&k_{22} &k_{23}\\
		0&0 &k_{33}\\
	\end{matrix}\right)\left(\begin{matrix}
		1&0 &0\\
		z^{\prime}_{21}&1&0\\
		z^{\prime}_{31}&z^{\prime}_{32}&1\\
	\end{matrix}\,\,\right)=:k_g\cdot (z\bar g).
\end{equation}
Here
\begin{equation}\label{d1}
	k_{11}=\frac{1}{\begin{vmatrix}
			g_{22}&g_{23}\\
			g_{32}&g_{33}\\
		\end{vmatrix}+\begin{vmatrix}
			g_{12}&g_{13}\\
			g_{32}&g_{33}\\
		\end{vmatrix}\cdot z_{21}+\begin{vmatrix}
			g_{12}&g_{13}\\
			g_{22}&g_{23}\\
		\end{vmatrix}\begin{vmatrix}
			z_{21}&1\\
			z_{31}&z_{32}\\
	\end{vmatrix}},\,\,\,\,\,\,\,\,
\end{equation}
\begin{equation}
	k_{22}=\frac{\begin{vmatrix}
			g_{22}&g_{23}\\
			g_{32}&g_{33}\\
		\end{vmatrix}+\begin{vmatrix}
			g_{12}&g_{13}\\
			g_{32}&g_{33}\\
		\end{vmatrix}\cdot z_{21}+\begin{vmatrix}
			g_{12}&g_{13}\\
			g_{22}&g_{23}\\
		\end{vmatrix}\begin{vmatrix}
			z_{21}&1\\
			z_{31}&z_{32}\\
	\end{vmatrix}}{g_{33}+g_{23}z_{32}+g_{13}z_{31}},\,\,\,\,\,\,\,\,\end{equation}
\begin{equation}\label{d3}
	k_{33}={g_{33}+g_{23}z_{32}+g_{13}z_{31}},\,\,\,\,\,\,\,\,\,\,\,\,\,\,\,\,\,\,\,\,\,\,\,\,\,\,\,\,\,\,\,\,\,\,\,\,\,\,\,\,\,\,\,\,\,\,\,\,\,\,\,\,\,\,\,\,\,\,\,\,\,\,\,\,\,\,\,\,\,\,\,\,\,\,\,\,\,\,\,\,
\end{equation}
\begin{equation}\label{z2}
	z^{\prime}_{21}=\frac{\begin{vmatrix}
			g_{21}&g_{23}\\
			g_{31}&g_{33}\\
		\end{vmatrix}+\begin{vmatrix}
			g_{11}&g_{13}\\
			g_{31}&g_{33}\\
		\end{vmatrix}\cdot z_{21}+\begin{vmatrix}
			g_{11}&g_{13}\\
			g_{21}&g_{23}\\
		\end{vmatrix}\begin{vmatrix}
			z_{21}&1\\
			z_{31}&z_{32}\\
	\end{vmatrix}}{\begin{vmatrix}
			g_{22}&g_{23}\\
			g_{32}&g_{33}\\
		\end{vmatrix}+\begin{vmatrix}
			g_{12}&g_{13}\\
			g_{32}&g_{33}\\
		\end{vmatrix}\cdot z_{21}+\begin{vmatrix}
			g_{12}&g_{13}\\
			g_{22}&g_{23}\\
		\end{vmatrix}\begin{vmatrix}
			z_{21}&1\\
			z_{31}&z_{32}\\
	\end{vmatrix}},\,\,\,\,\,\,\,\,
\end{equation}
\begin{equation}\label{z3}
	z^{\prime}_{31}=\frac{g_{31}+g_{21}z_{32}+g_{11}z_{31}}{g_{33}+g_{23}z_{32}+g_{13}z_{31}},\,\,\,\,\,\,\,\,z^{\prime}_{32}=\frac{g_{32}+g_{22}z_{32}+g_{12}z_{31}}{g_{33}+g_{23}z_{32}+g_{13}z_{31}}.\,
\end{equation}

Under the above notation, we may realize the principal series $T^{\chi}$ as follows. For $m_2,m_3\in\mathbb Z$, and $\rho_2,\rho_3\in\mathbb R$,  define a character $\chi=(m_2,m_3,\rho_2,\rho_3)$ of $K$ by
\begin{equation}
	\chi(k):=|k_{22}|^{m_2+\sqrt{-1}\rho_2}\cdot k_{22}^{-m_2}\cdot|k_{33}|^{m_3+\sqrt{-1}\rho_3}\cdot k_{33}^{-m_3}, \,\,k\in K.
\end{equation}
Denote by $\mathcal H_{\chi}$ the Hilbert space of square integrable functions on $Z$ with respect to the norm  
\begin{equation}
	(f,f):=\left(\frac{i}{2}\right)^3\int_Z |f(z)|^2\,dz_{21}\wedge d\overline{z_{21}} \wedge dz_{31}\wedge d\overline{z_{31}}\wedge dz_{32}\wedge d\overline{z_{32}}.%,\,\,f,h\in\mathcal H_{\chi}.
\end{equation}
Exploiting  (\ref{zzp})-(\ref{z3}), we can define the linear operator  $T^{\chi}(g)$ by
\begin{equation}
\begin{split}
&\left(T^{\chi}(g)f\right)(z):=\chi(k_g)\cdot\beta^{-1/2}(k_g)\cdot f(z\bar g)\\
&\,\,\,\,\,\,\,=|k_{22}|^{-2+m_2+\sqrt{-1}\rho_2}\cdot k_{22}^{-m_2}\cdot|k_{33}|^{-4+m_3+\sqrt{-1}\rho_3}\cdot k_{33}^{-m_3}\cdot f(z\bar g),\,\,f\in\mathcal H_{\chi}.
\end{split}
\end{equation}
By Theorem \ref{glnc}, $T^{\chi}$ is  unitary and irreducible;  $T^{\chi}$, $T^{\chi^{\prime}}$ are equivalently if and only if $(m^{\prime}_2,m^{\prime}_3,\rho^{\prime}_2,\rho^{\prime}_3)=(m_2,m_3,\rho_2,\rho_3)$, $(-m_2,m_3-m_2,-\rho_2,\rho_3-\rho_2)$,  $(m_3-m_2,-m_2,\rho_3-\rho_2,-\rho_2)$, $(-m_3,m_2-m_3,-\rho_3,\rho_2-\rho_3)$, $(m_3,m_2,\rho_3,\rho_2)$, or $(m_2-m_3,-m_3,\rho_2-\rho_3,-\rho_3)$. 
 
\smallskip

Similarly, we have the following complementary series. For  $m\in\mathbb Z$, $\rho\in\mathbb R$, and $0<\sigma<1$, define a character $\chi=(m,\rho,\sigma)$ by 
\begin{equation}
	\chi(k):=|k_{22}|^{m+\sqrt{-1}\rho+\sigma}\cdot k_{22}^{-m}\cdot|k_{33}|^{m+\sqrt{-1}\rho-\sigma}\cdot k_{33}^{-m}, \,\,k\in K.
\end{equation}
Define two unipotent subgroups of $SL_3(\mathbb C)$
\begin{equation}
\dot{Z}:=\left\{\left.\dot{z}=\left(\begin{matrix}
	1&0 &0\\
	0&1&0\\
	0&z_{1}&1\\
\end{matrix}\,\,\right)\right|z_{1}\in\mathbb C\,\right\},\,\,{\rm and}\,\,\,	\widehat Z:=\left\{\left. \hat{z}= \left(\begin{matrix}
1&0 &0\\
z_{21}&1&0\\
z_{31}&0&1\\
\end{matrix}\,\,\right)\right|z_{21},z_{31}\in\mathbb C\,\right\}.
\end{equation}
Since $Z=\dot{Z}\ltimes\widehat{Z}$, we write $z=(\hat{z},z_1)$.  Denote by $\mathcal H_{\chi}$ the Hilbert space of square integrable functions on $Z$ with respect to the norm
\begin{equation}
(f,f):=\left(\frac{i}{2}\right)^4\int_{\mathbb C\times Z} \frac{f(\hat{z},z_1)\cdot\overline{f(\hat{z},w_1)}}{\left|z_1-w_1\right|^{2-2\sigma}}\,dw_{1}\wedge d\overline{w_{1}}\wedge dz_{1}\wedge d\overline{z_{1}}\wedge  dz_{21}\wedge d\overline{z_{21}}\wedge dz_{31}\wedge d\overline{z_{31}}.%,\,\,f,h\in\mathcal H_{\sigma}.
\end{equation}
We define the linear operator $T^{\chi}(g)$ by 
\begin{equation}
\begin{split}
&\left(T^{\chi}(g)f\right)(z):=\chi(k_g)\cdot\beta^{-1/2}(k_g)\cdot f(z\bar g)\\
&\,\,\,\,\,\,\,\,\,\,\,=|k_{22}|^{-2+m+\sqrt{-1}\rho+\sigma}\cdot k_{22}^{-m}\cdot|k_{33}|^{-4+m+\sqrt{-1}\rho-\sigma}\cdot k_{33}^{-m}\cdot f(z\bar g),\,\, f\in\mathcal H_{\chi}.
\end{split}
\end{equation}
For different $\chi$, $T^{\chi}$ are unitary, irreducible, and inequivalent.

\begin{remark}
One can realize the  representation ${\rm Ind}_K^{SL_3(\mathbb C)}\chi$ by sections of the line bundle $V_{\chi}$ over  $SL_3(\mathbb C)/K$, where $V_{\chi}$ is the quotient space $\left(SL_3(\mathbb C)\times\mathbb C\right)/K$ by the equivalence relation
$(g, v)\sim \left(kg, \chi(k)\beta^{-1/2}(k)v\right)$,
$g\in SL_3(\mathbb C)$, $v\in\mathbb C$, $k\in K$. Notice that
\begin{equation}\label{smooths}
\mathcal D_{\chi}:=\left\{\begin{matrix}
{\rm smooth\,\,sections \,\,}\\
s: SL_3(\mathbb C)/K \rightarrow V_{\chi}
\end{matrix}\right\}\cong\left\{\begin{matrix}{\rm smooth\,\,sections \,\,}\\
\phi:SL_3(\mathbb C)\rightarrow \mathbb C
\end{matrix}\,\,\left|\begin{matrix}\,
\phi(kg)=\chi(k)\beta^{-1/2}(k)\phi(g)\\
\,{\rm for\,\,}g\in SL_3(\mathbb C)\,{\rm and}\,\, k \in K
\end{matrix} \right.\right\}.
\end{equation}
It is clear that $T^{\chi}(g)$ acts  on the right hand side of (\ref{smooths}) by	$\left(T^{\chi}(g)\phi\right)(z)=\phi(zg)$. Moreover, the Hilbert space $\mathcal H_{\chi}$ is the completion of $\mathcal D_{\chi}$, for both principal and complementary series, with respect to the metric on $V_{\chi}$ induced from the standard fiber metric $\frac{i}{2}dzd\bar z$ on $SL_3(\mathbb C)\times\mathbb C$. 
\end{remark}

At last, we introduce the degenerate principal series representation. Define
\begin{equation}
K_{2,1}:=\left\{\left.k=\left(\begin{matrix}
k_{11}&	k_{12}&	k_{13}\\
k_{21}&k_{22}&k_{23}\\
0&0&k_{33}\\
\end{matrix}\,\right)\right|k_{ij}\in\mathbb C,\,\,k_{33}\cdot \begin{vmatrix}
k_{11}&k_{12}\\
k_{21}&k_{22}\\
\end{vmatrix}=1\right\}.
\end{equation}
Equip  $K_{2,1}$ with a left Haar measure
\begin{equation}
			\small
d\mu_l(k):=\left(\frac{i}{2}\right)^{6}|k_{33}|^{6}   dk_{11}\wedge d\overline{k_{11}}\wedge dk_{12}\wedge d\overline{k_{12}}\wedge dk_{13}\wedge d\overline{k_{13}}\wedge dk_{21}\wedge d\overline{k_{21}}\wedge dk_{22}\wedge d\overline{k_{22}}\wedge dk_{23}\wedge d\overline{k_{23}},
\end{equation}
and a right Haar measure
\begin{equation}
d\mu_r(k):=\left(\frac{i}{2}\right)^{6} dk_{11}\wedge d\overline{k_{11}}\wedge dk_{12}\wedge d\overline{k_{12}}\wedge dk_{13}\wedge d\overline{k_{13}}\wedge dk_{21}\wedge d\overline{k_{21}}\wedge dk_{22}\wedge d\overline{k_{22}}\wedge dk_{23}\wedge d\overline{k_{23}}.
\end{equation}
Then the modular function $\beta$ is given by $\beta(k)=\frac{d\mu_l(k)}{d\mu_r(k)}=|k_{33}|^{6}$. Define
\begin{equation}
	Z_{2,1}:=\left\{\left.z=\left(\begin{matrix}
	1&0 &0\\
	0&1&0\\
	z_{31}&z_{32}&1\\
	\end{matrix}\,\,\right)\right|z_{31},z_{32}\in\mathbb C\,\right\}.
\end{equation}
Equip $Z_{2,1}$ with a Haar measure $d\mu(z)=(i/2)^2dz_{31}\wedge d\overline{z_{31}}\wedge dz_{32}\wedge d\overline{z_{32}}$.  It is clear that $Z_{2,1}$ is an open subset of the generalized flag variety $SL_3(\mathbb C)/K_{2,1}$.

For $m_2\in\mathbb Z$ and $\rho_2\in\mathbb R$, define a character $\chi=(m_2,\rho_2)$ of $K_{2,1}$ by
\begin{equation}
\chi(k):=|k_{33}|^{m_2+\sqrt{-1}\rho_2}\cdot k_{33}^{-m_2},\,\,\,\,k\in K_{2,1}.
\end{equation} 
Denote by $\mathcal H_{\chi}$  the Hilbert space of square integrable functions on $Z_{2,1}$ with respect to the norm
\begin{equation}
	(f,f):=\left(\frac{i}{2}\right)^2\int_{Z_{2,1}} |f(z)|^2\,dz_{31}\wedge d\overline{z_{31}}\wedge dz_{32}\wedge d\overline{z_{32}}.%,\,\,\,\,f,h\in\mathcal H_{\chi}.
\end{equation}
The linear operator  $T^{\chi}(g)$ acts on $\mathcal H_{\chi}$ by
\begin{equation}
\left(T^{\chi}(g)f\right)(z):=\chi(k_g)\cdot\beta^{-1/2}(k_g)\cdot f(z\bar g)=|k_{33}|^{-3+m_2+\sqrt{-1}\rho_2}\cdot k_{33}^{-m_2}\cdot f(z\bar g),\,\,f\in\mathcal H_{\chi}.
\end{equation}
Here $z\cdot g=k_g\cdot (z\bar g)$, for  generic $z\in Z_{2,1}$. More explicitly, 
\begin{equation}\label{zp}
	\left(\begin{matrix}
		1&0 &0\\
		0&1&0\\
		z_{31}&z_{32}&1\\
	\end{matrix}\,\,\right) \left(\begin{matrix}
		g_{11}&g_{12}&g_{13}\\
		g_{21}&g_{22}&g_{23}\\
		g_{31}&g_{32} &g_{33}\\
	\end{matrix}
	\right)=\left(\begin{matrix}
		k_{11}&	k_{12}&	k_{13}\\
		k_{21}&k_{22}&k_{23}\\
		0&0&k_{33}\\
	\end{matrix}\right)\left(\begin{matrix}
		1&0 &0\\
		0&1&0\\
		z^{\prime}_{31}&z^{\prime}_{32}&1\\
	\end{matrix}\,\,\right)=:k_g\cdot(z\bar g),
\end{equation}
where
\begin{equation}
	k_{11}=\frac{\begin{vmatrix}
			g_{11}&g_{13}\\
			g_{31}&g_{33}\\
		\end{vmatrix}+\begin{vmatrix}
			g_{11}&g_{13}\\
			g_{21}&g_{23}\\
		\end{vmatrix}\cdot z_{32}}{g_{33}+g_{23}z_{32}+g_{13}z_{31}},\,\,\,\,\,\,\,\,\,\,k_{12}=\frac{\begin{vmatrix}
		g_{12}&g_{13}\\
		g_{32}&g_{33}\\
	\end{vmatrix}+\begin{vmatrix}
	g_{12}&g_{13}\\
	g_{22}&g_{23}\\
\end{vmatrix}\cdot z_{32}}{g_{33}+g_{23}z_{32}+g_{13}z_{31}},\,
\end{equation}
\begin{equation}\label{k2}
	k_{21}=\frac{\begin{vmatrix}
			g_{21}&g_{23}\\
			g_{31}&g_{33}\\
		\end{vmatrix}-\begin{vmatrix}
			g_{11}&g_{13}\\
			g_{21}&g_{23}\\
		\end{vmatrix}\cdot z_{31}}{g_{33}+g_{23}z_{32}+g_{13}z_{31}},\,\,\,\,\,\,\,\,\,\,k_{22}=\frac{\begin{vmatrix}
			g_{22}&g_{23}\\
			g_{32}&g_{33}\\
		\end{vmatrix}-\begin{vmatrix}
			g_{12}&g_{13}\\
			g_{22}&g_{23}\\
		\end{vmatrix}\cdot z_{31}}{g_{33}+g_{23}z_{32}+g_{13}z_{31}},
\end{equation}
\begin{equation}
	k_{13}=g_{13},\,\,\,\,\,\,\,\,\,\,\,\,\,\,\,\,k_{23}=g_{23},\,\,\,\,\,\,\,\,\,\,\,\,\,\,\,\,k_{33}=g_{33}+g_{23}z_{32}+g_{13}z_{31},\,\,\,\,\,\,\,\,\,\,\,\,\,\,\,\,\,\,\,\,\,\,\,\,\,\,\,\,\,
\end{equation}
\begin{equation}\label{dd1}
	z_{31}^{\prime}=\frac{g_{31}+g_{21}z_{32}+g_{11}z_{31}}{g_{33}+g_{23}z_{32}+g_{13}z_{31}},\,\,\,\,\,\,\,\,\,\,\,\,z_{32}^{\prime}=\frac{g_{32}+g_{22}z_{32}+g_{12}z_{31}}{g_{33}+g_{23}z_{32}+g_{13}z_{31}}.\,\,\,\,\,\,\,\,\,\,\,\,\,\,\,\,\,\,\,\,\,\,\,\,\,\,\,\,
\end{equation}
For different $\chi$, $T^{\chi}$ are unitary, irreducible, and inequivalent.

\section{Irreducible unitary representations of \texorpdfstring{$SL_3(\mathbb R)$} {rr}}
We denote by $\mathring K$ the upper triangular matrix subgroup of $SL_3(\mathbb R)$ consisting of matrices $x=(x_{pq})$ where $x_{pq}=0$ for $p>q$. Denote by $X$ the unipotent lower triangular subgroup of $SL_3(\mathbb R)$ consisting of matrices $x=(x_{pq})$ where $x_{pp}=1$ and $x_{pq}=0$ for $p<q$. Equip $\mathring K$ with a left Haar measure
\begin{equation}
	d\mu_l(k):=|k_{33}|\cdot dk_{12}\wedge dk_{13}\wedge dk_{22}\wedge dk_{23}\wedge dk_{33},
\end{equation}
and a right Haar measure 
\begin{equation}
	d\mu_r(k):=|k_{22}|^{-2}\cdot|k_{33}|^{-3}\cdot dk_{12}\wedge dk_{13}\wedge dk_{22}\wedge dk_{23}\wedge dk_{33}.
\end{equation}
The modular function of $\mathring K$ is
\begin{equation}
	\beta(k):=\frac{d\mu_l(k)}{d\mu_r(k)}=|k_{22}|^{2}\cdot |k_{33}|^{4}.
\end{equation}
Each $g\in SL_3(\mathbb R)$ induces a birational transformation $\bar g:x\dashmapsto x\bar g$ of $X$, where $x\bar g$ are defined by the decomposition $x\cdot g=k_g\cdot(x\bar g)$. Notice that the explicit formulas for $x\bar g\in X$ and   $k_g\in \mathring K$ can be derived by substituting   $x=(x_{pq})$ for $z=(z_{pq})$ in (\ref{zzp})-(\ref{z3}). For $m_2,m_3\in\{0,1\}$, and $\rho_2,\rho_3\in\mathbb R$, define a character  $\chi=(m_2,m_3,\rho_2,\rho_3)$ of $\mathring K$ by
\begin{equation}
\chi(k):=|k_{22}|^{\sqrt{-1} \rho_2} \cdot|k_{33}|^{\sqrt{-1} \rho_3} \cdot\left(\frac{k_{22}}{|k_{22}|}\right)^{m_2}\cdot\left(\frac{k_{33}}{|k_{33}|}\right)^{m_3},\,\,k\in \mathring K.
\end{equation}
Denote by $\mathcal H_{\chi}$ the Hilbert space of square integrable functions on $X$ with respect to the norm  
\begin{equation}
	(f,f):=\int_X |f(x)|^2\,dx_{21}\,  dx_{31}\,dx_{32}.%,\,\,f,h\in\mathcal H_{\chi}.
\end{equation}
Let
\begin{equation}
	\begin{split}
		&\left(T^{\chi}(g)f\right)(x):=\chi(k_g)\cdot\beta^{-1/2}(k_g)\cdot f(x\bar g)\\
		&\,\,\,\,\,\,\,=|k_{22}|^{-1+\sqrt{-1} \rho_2} \cdot|k_{33}|^{-2+\sqrt{-1} \rho_3} \cdot\left(\frac{k_{22}}{|k_{22}|}\right)^{m_2}\cdot\left(\frac{k_{33}}{|k_{33}|}\right)^{m_3}\cdot f(x\bar g),\,\,f\in\mathcal  H_{\chi}.
	\end{split}
\end{equation}
Then by Theorem \ref{glnr}, $T^{\chi}$ is  unitary and irreducible;  $T^{\chi}$, $T^{\chi^{\prime}}$ are equivalently if and only if $(m^{\prime}_2,m^{\prime}_3,\rho^{\prime}_2,\rho^{\prime}_3)=(m_2,m_3,\rho_2,\rho_3)$, $(-m_2,m_3-m_2,-\rho_2,\rho_3-\rho_2)$,  $(m_3-m_2,-m_2,\rho_3-\rho_2,-\rho_2)$, $(-m_3,m_2-m_3,-\rho_3,\rho_2-\rho_3)$, $(m_3,m_2,\rho_3,\rho_2)$, or $(m_2-m_3,-m_3,\rho_2-\rho_3,-\rho_3)$. 

\smallskip

Next, we have the following complementary series. For $m\in\{0,1\}$, $\rho\in\mathbb R$, and $0<\sigma<1$, define a character $\chi=(m,\rho,\sigma)$ by 
\begin{equation}
	\chi(k):=|k_{22}|^{\sqrt{-1}\rho+\sigma}\cdot \left(\frac{k_{22}}{|k_{22}|}\right)^{m}\cdot|k_{33}|^{\sqrt{-1}\rho-\sigma}\cdot \left(\frac{k_{33}}{|k_{33}|}\right)^{m}, \,\,k\in K.
\end{equation}
Let
\begin{equation}
	\dot{X}:=\left\{\left.\dot{x}=\left(\begin{matrix}
	1&0 &0\\
	0&1&0\\
	0&x_{1}&1\\
\end{matrix}\,\,\right)\right|x_{1}\in\mathbb R\,\right\}\,\,{\rm and}\,\,\,\widehat X:=\left\{\left. \hat{x}= \left(\begin{matrix}
1&0 &0\\
x_{21}&1&0\\
x_{31}&0&1\\
\end{matrix}\,\,\right)\right|x_{21},x_{31}\in\mathbb R\,\right\}.
\end{equation}
Since $X=\dot{X}\ltimes\widehat{X}$, we may write $x=(\hat{x},x_1)$.  Denote by $\mathcal H_{\chi}$ the Hilbert space of square integrable functions on $X$ with respect to the norm
\begin{equation}
(f,f):=\int_{\mathbb R\times X} \frac{f(\hat{x},x_1)\cdot\overline{f(\hat{x},y_1)}}{\left|x_1-y_1\right|^{1-\sigma}}\,dx_{1}\, dy_{1}\,  dx_{21}\, dx_{31}.%,\,\,f,h\in\mathcal H_{\sigma}.
\end{equation}
Define 
\begin{equation}
\left(T^{\chi}(g)f\right)(z):=|k_{22}|^{-1+\sqrt{-1}\rho+\sigma}\cdot\left(\frac{k_{22}}{|k_{22}|}\right)^{-m}\cdot|k_{33}|^{-2+\sqrt{-1}\rho-\sigma}\cdot\left(\frac{k_{33}}{|k_{33}|}\right)^{-m}\cdot f(x\bar g),\,\, f\in\mathcal H_{\chi}.
\end{equation}
For different $\chi$, $T^{\chi}$ are unitary, irreducible, and inequivalent.

\smallskip

In a parallel way, we may introduce the following degenerate principal series. Define
\begin{equation}
	\mathring K_{2,1}:=\left\{\left.k=\left(\begin{matrix}
		k_{11}&	k_{12}&	k_{13}\\
		k_{21}&k_{22}&k_{23}\\
		0&0&k_{33}\\
	\end{matrix}\,\right)\right|k_{ij}\in\mathbb R,\,\,k_{33}\cdot \begin{vmatrix}
		k_{11}&k_{12}\\
		k_{21}&k_{22}\\
	\end{vmatrix}=1\right\}.
\end{equation}
Equip  $\mathring K_{2,1}$ with a left Haar measure
\begin{equation}
	d\mu_l(k):= |k_{33}|^{3}   dk_{11}\wedge dk_{12}\wedge  dk_{13}\wedge  dk_{21}\wedge  dk_{22}\wedge  dk_{23},
\end{equation}
and a right Haar measure
\begin{equation}
	d\mu_r(k):=dk_{11}\wedge  dk_{12}\wedge  dk_{13}\wedge  dk_{21}\wedge  dk_{22}\wedge  dk_{23}.
\end{equation}
Then the modular function $\beta(k)=|k_{33}|^{3}$. Define
\begin{equation}
	X_{2,1}:=\left\{\left.x=\left(\begin{matrix}
		1&0 &0\\
		0&1&0\\
		x_{31}&x_{32}&1\\
	\end{matrix}\,\,\right)\right|x_{31},x_{32}\in\mathbb R\,\right\}.
\end{equation}
For  $m_2\in\{0,1\}$ and $\rho_2\in\mathbb R$, define a character $\chi=(m_2,\rho_2)$ of $K_{2,1}$ by
\begin{equation}
	\chi(k):=|k_{33}|^{m_2+\sqrt{-1}\rho_2}\cdot k_{33}^{-m_2},\,\,\,\,k\in K_{2,1}.
\end{equation} 
Denote by $\mathcal H_{\chi}$  the Hilbert space of square integrable functions on $X_{2,1}$ with respect to the norm
\begin{equation}
	(f,f):=\int_{X_{2,1}} |f(x)|^2\,dx_{31}\wedge dx_{32}.%,\,\,\,\,f,h\in\mathcal H_{\chi}.
\end{equation}
Exploiting the decomposition $x\cdot g=k_g\cdot (x\bar g)$ for  generic $x\in X_{2,1}$, we may define
\begin{equation}
	\left(T^{\chi}(g)f\right)(x):=\chi(k_g)\cdot\beta^{-1/2}(k_g)\cdot f(x\bar g)=|k_{33}|^{-\frac{3}{2}+m_2+\sqrt{-1}\rho_2}\cdot k_{33}^{-m_2}\cdot f(x\bar g),\,\,f\in\mathcal H_{\chi}.
\end{equation}
Notice that the explicit formulas for $k_g\in\mathring K_{2,1}$ and $x\bar g\in X_{2,1}$ can be given by substituting $x=(x_{ij})$ for $z=(z_{ij})$ in (\ref{zp})-(\ref{dd1}). For different $\chi$, $T^{\chi}$ are unitary, irreducible, and inequivalent. 

Thanks to the appearance of  non-conjugate maximal Abelian subgroups of $SL_3(\mathbb R)$,  we have one further principal series as shown by \cite{GG1}. Define subsets
\begin{equation}
	\dot Z_{1,1}:=\left\{\left.\dot{z}=\left(\begin{matrix}
		1&0 &0\\
		z_{21}&1&0\\
		0&0&1\\
	\end{matrix}\right)\right|\begin{matrix}
	z_{21}\in\mathbb C\\
	\Im z_{21}\neq 0
\end{matrix}\right\},\,\,{X}_{1,1}:=\left\{\left. {x}= \left(\begin{matrix}
		1&0 &0\\
		0&1&0\\
		x_{31}&x_{32}&1\\
	\end{matrix}\,\right)\right|x_{31},x_{32}\in\mathbb R\,\right\},
\end{equation}
and 
\begin{equation}
Z_{1}:=\dot Z_{1,1}\cdot X_{1,1}=\left\{\dot{z}\cdot x\left|\dot{z}\in\dot{Z}_{1,1},\,x\in X_{1,1}\right.\right\}.
\end{equation}
It is clear that $Z_{1}$ is a subset of $Z$ and hence can be viewed as a submanifold of $SL_3(\mathbb C)/K$ (see \S \ref{sl3c} for the definitions of $Z$ and $K$). One can further verify that $Z_1$ is transitive with respect to the birational action of $SL_3(\mathbb R)$.  Write $\chi=(n_1,\rho_1)$ for $n_1\in\mathbb Z^{+}$ and $\rho_1\in\mathbb R$. Denote by $\mathcal H_{\chi}$ the Hilbert space of functions $f(z_{21},x_{31},x_{32})$ defined on $Z_1$, which are analytic in the variable $z_{21}=x_{21}+\sqrt{-1}y_{21}$ on the upper half-plane $\left\{\Im z_{21}>0\right\}$ and the lower planes $\left\{\Im z_{21}<0\right\}$, and are square integrable with respect to the following norm.  
\begin{equation}\label{n11}
	(f,f):=\int_{Z_1} |f(z_{21},x_{31},x_{32})|^2\cdot \frac{\left|\Im z_{21}\right|^{n_1-2}}{\Gamma(n_1-1)}\,dx_{21}\, dy_{21}\,dx_{31}\,dx_{32}.%,\,\,f,h\in\mathcal H_{\chi}.
\end{equation}
Define the linear operator  $T^{\chi}(g)$ acting on $f\in\mathcal  H_{\chi}$ by
\begin{equation}\label{tn11}
	\small
\left(T^{\chi}(g)f\right)(z_{21},x_{31},x_{32}):=\left|\alpha_1\delta_1-\beta_1\gamma_1\right|^{\frac{3}{2}+\frac{n_1}{2}+\sqrt{-1}\rho_1}\left(\beta_1z_{21}+\delta_1\right)^{-n_1} f\left(\frac{\alpha_1z_{21}+\gamma_1}{\beta_1z_{21}+\delta_1},x_{31},x_{32}\right),
\end{equation}
where $\alpha_1,\beta_1,\gamma_1,\delta_1$ are determined by (\ref{zp})-(\ref{k2})as follows. 
\begin{equation}
	\left(\begin{matrix}
		1&0 &0\\
		0&1&0\\
		x_{31}&x_{32}&1\\
	\end{matrix}\,\,\right) \left(\begin{matrix}
		g_{11}&g_{12}&g_{13}\\
		g_{21}&g_{22}&g_{23}\\
		g_{31}&g_{32} &g_{33}\\
	\end{matrix}
	\right)=\left(\begin{matrix}
		\alpha_1&\beta_1&	k_{13}\\
		\gamma_1&\delta_1&k_{23}\\
		0&0&k_{33}\\
	\end{matrix}\right)\left(\begin{matrix}
		1&0 &0\\
		0&1&0\\
		x^{\prime}_{31}&x^{\prime}_{32}&1\\
	\end{matrix}\,\,\right),
\end{equation}
where
\begin{equation}
	\alpha_{1}=\frac{\begin{vmatrix}
			g_{11}&g_{13}\\
			g_{31}&g_{33}\\
		\end{vmatrix}+\begin{vmatrix}
			g_{11}&g_{13}\\
			g_{21}&g_{23}\\
		\end{vmatrix}\cdot x_{32}}{g_{33}+g_{23}x_{32}+g_{13}x_{31}},\,\,\,\,\,\,\,\,\,\,\beta_{1}=\frac{\begin{vmatrix}
			g_{12}&g_{13}\\
			g_{32}&g_{33}\\
		\end{vmatrix}+\begin{vmatrix}
			g_{12}&g_{13}\\
			g_{22}&g_{23}\\
		\end{vmatrix}\cdot x_{32}}{g_{33}+g_{23}x_{32}+g_{13}x_{31}},\,
\end{equation}
\begin{equation}
	\gamma_{1}=\frac{\begin{vmatrix}
			g_{21}&g_{23}\\
			g_{31}&g_{33}\\
		\end{vmatrix}-\begin{vmatrix}
			g_{11}&g_{13}\\
			g_{21}&g_{23}\\
		\end{vmatrix}\cdot x_{31}}{g_{33}+g_{23}x_{32}+g_{13}x_{31}},\,\,\,\,\,\,\,\,\,\,\delta_{1}=\frac{\begin{vmatrix}
			g_{22}&g_{23}\\
			g_{32}&g_{33}\\
		\end{vmatrix}-\begin{vmatrix}
			g_{12}&g_{13}\\
			g_{22}&g_{23}\\
		\end{vmatrix}\cdot x_{31}}{g_{33}+g_{23}x_{32}+g_{13}x_{31}},
\end{equation}
and hence (see (\ref{z2}) for reference)
\begin{equation}
\frac{\alpha_1z_{21}+\gamma_1}{\beta_1z_{21}+\delta_1}=\frac{\begin{vmatrix}
		g_{21}&g_{23}\\
		g_{31}&g_{33}\\
	\end{vmatrix}+\begin{vmatrix}
		g_{11}&g_{13}\\
		g_{31}&g_{33}\\
	\end{vmatrix}\cdot z_{21}+\begin{vmatrix}
		g_{11}&g_{13}\\
		g_{21}&g_{23}\\
	\end{vmatrix}\begin{vmatrix}
		z_{21}&1\\
		x_{31}&x_{32}\\
\end{vmatrix}}{\begin{vmatrix}
		g_{22}&g_{23}\\
		g_{32}&g_{33}\\
	\end{vmatrix}+\begin{vmatrix}
		g_{12}&g_{13}\\
		g_{32}&g_{33}\\
	\end{vmatrix}\cdot z_{21}+\begin{vmatrix}
		g_{12}&g_{13}\\
		g_{22}&g_{23}\\
	\end{vmatrix}\begin{vmatrix}
		z_{21}&1\\
		x_{31}&x_{32}\\
\end{vmatrix}}.\,\,\,\,\,\,\,\,
\end{equation}
For different $\chi$, $T^{\chi}$ are  unitary, irreducible, and inequivalent.

\section{Irreducible unitary representations of \texorpdfstring{$SL(4,\mathbb C)$} {rr}}

Similarly, we denote by $K$ the upper triangular matrix subgroup of $SL_4(\mathbb C)$, which consists of complex matrices $k=(k_{pq})$ such that $k_{pq}=0$ for $p>q$. Denote by $Z$ the unipotent lower triangular matrix subgroup of $SL_4(\mathbb C)$, which consists of matrices $z=(z_{pq})$ such that $z_{pp}=1$ and $z_{pq}=0$ for $p<q$. %Equip $K$ with a left Haar measure\begin{equation}d\mu_l(k):=\left(\frac{i}{2}\right)^{5}|k_{33}|^2\cdot dk_{12}\wedge d\overline{k_{12}}\wedge dk_{13}\wedge d\overline{k_{13}}\wedge dk_{22}\wedge d\overline{k_{22}}\wedge dk_{23}\wedge d\overline{k_{23}}\wedge dk_{33}\wedge d\overline{k_{33}},\end{equation}and a right Haar measure \begin{equation}d\mu_r(k):=\left(\frac{i}{2}\right)^{5}|k_{22}|^{-4}\cdot|k_{33}|^{-6}\cdot dk_{12}\wedge d\overline{k_{12}}\wedge dk_{13}\wedge d\overline{k_{13}}\wedge dk_{22}\wedge d\overline{k_{22}}\wedge dk_{23}\wedge d\overline{k_{23}}\wedge dk_{33}\wedge d\overline{k_{33}}.\end{equation}
Then the modular function of  $K$ is  
\begin{equation}\label{beta2}
	\beta(k)=\frac{d\mu_l(k)}{d\mu_r(k)}=|k_{22}|^{4}\cdot |k_{33}|^{8}\cdot |k_{44}|^{12}.
\end{equation}
For $g\in SL_3(\mathbb C)$ and generic $z\in Z$, we have
\begin{equation}\label{zzzp}
	\footnotesize
	\left(\begin{matrix}
		1&0 &0&0\\
		z_{21}&1&0&0\\
		z_{31}&z_{32}&1&0\\
		z_{41}&z_{42}&z_{43}&1\\
	\end{matrix}\,\right) \left(\begin{matrix}
		g_{11}&g_{12}&g_{13}&g_{14}\\
		g_{21}&g_{22}&g_{23}&g_{24}\\
		g_{31}&g_{32} &g_{33}&g_{34}\\
		g_{41}&g_{42} &g_{43}&g_{44}\\
	\end{matrix}
	\right)=\left(\begin{matrix}k_{11}&k_{12} &k_{13}&k_{14} \\0&k_{22} &k_{23}&k_{24}\\0&0 &k_{33}&k_{34}\\
		0&0 &0&k_{44}\\
	\end{matrix}\right)\left(\begin{matrix}
		1&0 &0&0\\
		z^{\prime}_{21}&1&0&0\\
		z^{\prime}_{31}&z^{\prime}_{32}&1&0\\
		z^{\prime}_{41}&z^{\prime}_{42}&z^{\prime}_{43}&1\\
	\end{matrix}\,\right)=:k_g\cdot z\bar g\,,
\end{equation}
where
\begin{equation}
 	\footnotesize
 	z^{\prime}_{21}=\frac{\begin{vmatrix} g_{12}&g_{13}&g_{14}\\			g_{32}&g_{33}&g_{34}\\
			g_{42}&g_{43}&g_{44}\\
		\end{vmatrix}+\begin{vmatrix} g_{11}&g_{13}&g_{14}\\
			g_{31}&g_{33}&g_{34}\\
			g_{41}&g_{43}&g_{44}\\
		\end{vmatrix}z_{21}+\begin{vmatrix} g_{11}&g_{12}&g_{14}\\
			g_{31}&g_{32}&g_{34}\\
			g_{41}&g_{42}&g_{44}\\
		\end{vmatrix}\begin{vmatrix} z_{21}&1\\
			z_{31}&z_{32}\\
		\end{vmatrix}+\begin{vmatrix} g_{11}&g_{12}&g_{13}\\
			g_{31}&g_{32}&g_{33}\\
			g_{41}&g_{42}&g_{43}\\
		\end{vmatrix}\begin{vmatrix} z_{21}&1&0\\
			z_{31}&z_{32}&1\\
			z_{41}&z_{42}&z_{43}\\
	\end{vmatrix}}{\begin{vmatrix} g_{22}&g_{23}&g_{24}\\
			g_{32}&g_{33}&g_{34}\\
			g_{42}&g_{43}&g_{44}\\
		\end{vmatrix}+\begin{vmatrix} g_{21}&g_{23}&g_{24}\\
			g_{31}&g_{33}&g_{34}\\
			g_{41}&g_{43}&g_{44}\\
		\end{vmatrix}z_{21}+\begin{vmatrix} g_{21}&g_{22}&g_{24}\\
			g_{31}&g_{32}&g_{34}\\
			g_{41}&g_{42}&g_{44}\\
		\end{vmatrix}\begin{vmatrix} z_{21}&1\\
			z_{31}&z_{32}\\
		\end{vmatrix}+\begin{vmatrix} g_{21}&g_{22}&g_{23}\\
			g_{31}&g_{32}&g_{33}\\
			g_{41}&g_{42}&g_{43}\\
		\end{vmatrix}\begin{vmatrix} z_{21}&1&0\\
			z_{31}&z_{32}&1\\
			z_{41}&z_{42}&z_{43}\\
	\end{vmatrix}}\,\,,
\end{equation}
\begin{equation}\label{z31}
	\scriptsize
	z^{\prime}_{31}=\frac{\begin{vmatrix}
			g_{13}&g_{14}\\
			g_{43}&g_{44}\\
		\end{vmatrix}+\begin{vmatrix}
			g_{12}&g_{14}\\
			g_{42}&g_{44}\\
		\end{vmatrix} z_{32}+\begin{vmatrix}
			g_{11}&g_{14}\\
			g_{41}&g_{44}\\
		\end{vmatrix} z_{31}+\begin{vmatrix}
			g_{12}&g_{13}\\
			g_{42}&g_{43}\\
		\end{vmatrix}\begin{vmatrix}
			z_{32}&1\\
			z_{42}&z_{43}\\
		\end{vmatrix}+\begin{vmatrix}
			g_{11}&g_{13}\\
			g_{41}&g_{43}\\
		\end{vmatrix}\begin{vmatrix}
			z_{31}&1\\
			z_{41}&z_{43}\\
		\end{vmatrix}+\begin{vmatrix}
			g_{11}&g_{12}\\
			g_{41}&g_{42}\\
		\end{vmatrix}\begin{vmatrix}
			z_{31}&z_{32}\\
			z_{41}&z_{42}\\
	\end{vmatrix}}{\begin{vmatrix}
			g_{33}&g_{34}\\
			g_{43}&g_{44}\\
		\end{vmatrix}+\begin{vmatrix}
			g_{32}&g_{34}\\
			g_{42}&g_{44}\\
		\end{vmatrix} z_{32}+\begin{vmatrix}
			g_{31}&g_{34}\\
			g_{41}&g_{44}\\
		\end{vmatrix} z_{31}+\begin{vmatrix}
			g_{32}&g_{33}\\
			g_{42}&g_{43}\\
		\end{vmatrix}\begin{vmatrix}
			z_{32}&1\\
			z_{42}&z_{43}\\
		\end{vmatrix}+\begin{vmatrix}
			g_{31}&g_{33}\\
			g_{41}&g_{43}\\
		\end{vmatrix}\begin{vmatrix}
			z_{31}&1\\
			z_{41}&z_{43}\\
		\end{vmatrix}+\begin{vmatrix}
			g_{31}&g_{32}\\
			g_{41}&g_{42}\\
		\end{vmatrix}\begin{vmatrix}
			z_{31}&z_{32}\\
			z_{41}&z_{42}\\
	\end{vmatrix}}\,\,,
\end{equation}

\begin{equation}\label{z32}
	\scriptsize
	z^{\prime}_{32}=\frac{\begin{vmatrix}
			g_{23}&g_{24}\\
			g_{43}&g_{44}\\
		\end{vmatrix}+\begin{vmatrix}
			g_{22}&g_{24}\\
			g_{42}&g_{44}\\
		\end{vmatrix} z_{32}+\begin{vmatrix}
			g_{21}&g_{24}\\
			g_{41}&g_{44}\\
		\end{vmatrix} z_{31}+\begin{vmatrix}
			g_{22}&g_{23}\\
			g_{42}&g_{43}\\
		\end{vmatrix}\begin{vmatrix}
			z_{32}&1\\
			z_{42}&z_{43}\\
		\end{vmatrix}+\begin{vmatrix}
			g_{21}&g_{23}\\
			g_{41}&g_{43}\\
		\end{vmatrix}\begin{vmatrix}
			z_{31}&1\\
			z_{41}&z_{43}\\
		\end{vmatrix}+\begin{vmatrix}
			g_{21}&g_{22}\\
			g_{41}&g_{42}\\
		\end{vmatrix}\begin{vmatrix}
			z_{31}&z_{32}\\
			z_{41}&z_{42}\\
	\end{vmatrix}}{\begin{vmatrix}
			g_{33}&g_{34}\\
			g_{43}&g_{44}\\
		\end{vmatrix}+\begin{vmatrix}
			g_{32}&g_{34}\\
			g_{42}&g_{44}\\
		\end{vmatrix} z_{32}+\begin{vmatrix}
			g_{31}&g_{34}\\
			g_{41}&g_{44}\\
		\end{vmatrix} z_{31}+\begin{vmatrix}
			g_{32}&g_{33}\\
			g_{42}&g_{43}\\
		\end{vmatrix}\begin{vmatrix}
			z_{32}&1\\
			z_{42}&z_{43}\\
		\end{vmatrix}+\begin{vmatrix}
			g_{31}&g_{33}\\
			g_{41}&g_{43}\\
		\end{vmatrix}\begin{vmatrix}
			z_{31}&1\\
			z_{41}&z_{43}\\
		\end{vmatrix}+\begin{vmatrix}
			g_{31}&g_{32}\\
			g_{41}&g_{42}\\
		\end{vmatrix}\begin{vmatrix}
			z_{31}&z_{32}\\
			z_{41}&z_{42}\\
	\end{vmatrix}}\,\,,
\end{equation}

\begin{equation}\label{z41}
	\scriptsize
	z^{\prime}_{41}=\frac{g_{14}+g_{13}z_{43}+g_{12}z_{42}+g_{11}z_{41}}{g_{44}+g_{43}z_{43}+g_{42}z_{42}+g_{41}z_{41}},\,\,z^{\prime}_{42}=\frac{g_{24}+g_{23}z_{43}+g_{32}z_{42}+g_{31}z_{41}}{g_{44}+g_{43}z_{43}+g_{42}z_{42}+g_{41}z_{41}},\,\,z^{\prime}_{43}=\frac{g_{34}+g_{33}z_{43}+g_{32}z_{42}+g_{31}z_{41}}{g_{44}+g_{43}z_{43}+g_{42}z_{42}+g_{41}z_{41}}\,\,,
\end{equation}
and
\begin{equation}
	\scriptsize
	k_{11}=\frac{1}{\begin{vmatrix} g_{22}&g_{23}&g_{24}\\
			g_{32}&g_{33}&g_{34}\\
			g_{42}&g_{43}&g_{44}\\
		\end{vmatrix}+\begin{vmatrix} g_{21}&g_{23}&g_{24}\\
			g_{31}&g_{33}&g_{34}\\
			g_{41}&g_{43}&g_{44}\\
		\end{vmatrix}z_{21}+\begin{vmatrix} g_{21}&g_{22}&g_{24}\\
			g_{31}&g_{32}&g_{34}\\
			g_{41}&g_{42}&g_{44}\\
		\end{vmatrix}\begin{vmatrix} z_{21}&1\\
			z_{31}&z_{32}\\
		\end{vmatrix}+\begin{vmatrix} g_{21}&g_{22}&g_{23}\\
			g_{31}&g_{32}&g_{33}\\
			g_{41}&g_{42}&g_{43}\\
		\end{vmatrix}\begin{vmatrix} z_{21}&1&0\\
			z_{31}&z_{32}&1\\
			z_{41}&z_{42}&z_{43}\\
	\end{vmatrix}}\,\,,
\end{equation}
\begin{equation}
	\scriptsize
	k_{22}=\frac{\begin{vmatrix} g_{22}&g_{23}&g_{24}\\
			g_{32}&g_{33}&g_{34}\\
			g_{42}&g_{43}&g_{44}\\
		\end{vmatrix}+\begin{vmatrix} g_{21}&g_{23}&g_{24}\\
			g_{31}&g_{33}&g_{34}\\
			g_{41}&g_{43}&g_{44}\\
		\end{vmatrix}z_{21}+\begin{vmatrix} g_{21}&g_{22}&g_{24}\\
			g_{31}&g_{32}&g_{34}\\
			g_{41}&g_{42}&g_{44}\\
		\end{vmatrix}\begin{vmatrix} z_{21}&1\\
			z_{31}&z_{32}\\
		\end{vmatrix}+\begin{vmatrix} g_{21}&g_{22}&g_{23}\\
			g_{31}&g_{32}&g_{33}\\
			g_{41}&g_{42}&g_{43}\\
		\end{vmatrix}\begin{vmatrix} z_{21}&1&0\\
			z_{31}&z_{32}&1\\
			z_{41}&z_{42}&z_{43}\\
	\end{vmatrix}}{\begin{vmatrix}
			g_{33}&g_{34}\\
			g_{43}&g_{44}\\
		\end{vmatrix}+\begin{vmatrix}
			g_{32}&g_{34}\\
			g_{42}&g_{44}\\
		\end{vmatrix} z_{32}+\begin{vmatrix}
			g_{31}&g_{34}\\
			g_{41}&g_{44}\\
		\end{vmatrix} z_{31}+\begin{vmatrix}
			g_{32}&g_{33}\\
			g_{42}&g_{43}\\
		\end{vmatrix}\begin{vmatrix}
			z_{32}&1\\
			z_{42}&z_{43}\\
		\end{vmatrix}+\begin{vmatrix}
			g_{31}&g_{33}\\
			g_{41}&g_{43}\\
		\end{vmatrix}\begin{vmatrix}
			z_{31}&1\\
			z_{41}&z_{43}\\
		\end{vmatrix}+\begin{vmatrix}
			g_{31}&g_{32}\\
			g_{41}&g_{42}\\
		\end{vmatrix}\begin{vmatrix}
			z_{31}&z_{32}\\
			z_{41}&z_{42}\\
	\end{vmatrix}}\,\,,\,\,\end{equation}
\begin{equation}\label{k33}
	\scriptsize
	k_{33}=\frac{\begin{vmatrix}
			g_{33}&g_{34}\\
			g_{43}&g_{44}\\
		\end{vmatrix}+\begin{vmatrix}
			g_{32}&g_{34}\\
			g_{42}&g_{44}\\
		\end{vmatrix} z_{32}+\begin{vmatrix}
			g_{31}&g_{34}\\
			g_{41}&g_{44}\\
		\end{vmatrix} z_{31}+\begin{vmatrix}
			g_{32}&g_{33}\\
			g_{42}&g_{43}\\
		\end{vmatrix}\begin{vmatrix}
			z_{32}&1\\
			z_{42}&z_{43}\\
		\end{vmatrix}+\begin{vmatrix}
			g_{31}&g_{33}\\
			g_{41}&g_{43}\\
		\end{vmatrix}\begin{vmatrix}
			z_{31}&1\\
			z_{41}&z_{43}\\
		\end{vmatrix}+\begin{vmatrix}
			g_{31}&g_{32}\\
			g_{41}&g_{42}\\
		\end{vmatrix}\begin{vmatrix}
			z_{31}&z_{32}\\
			z_{41}&z_{42}\\
	\end{vmatrix}}{g_{44}+g_{43}z_{43}+g_{42}z_{42}+g_{41}z_{41}}\,,
\end{equation}
\begin{equation}\label{k44}
	\footnotesize
	k_{44}={g_{44}+g_{43}z_{43}+g_{42}z_{42}+g_{41}z_{41}}\,.\,\,\,\,\,\,\,\,\,\,\,\,\,\,\,\,\,\,\,\,\,\,\,\,\,\,\,\,\,\,\,\,\,\,\,\,\,\,\,\,\,\,\,\,\,\,\,\,\,\,\,\,\,\,\,\,\,\,\,\,\,\,\,\,\,\,\,\,\,\,\,\,\,\,\,\,\,\,\,\,
\end{equation}

We may realize the principal series $T^{\chi}$ as follows. For $m_2,m_3\in\mathbb Z$, and $\rho_2,\rho_3\in\mathbb R$,  define a character $\chi=(m_2,m_3,m_4,\rho_2,\rho_3,\rho_4)$ of $K$ by
\begin{equation}
	\chi(k):=|k_{22}|^{m_2+\sqrt{-1}\rho_2}\cdot k_{22}^{-m_2}\cdot|k_{33}|^{m_3+\sqrt{-1}\rho_3}\cdot k_{33}^{-m_3} \cdot|k_{44}|^{m_4+\sqrt{-1}\rho_4}\cdot k_{44}^{-m_4}, \,\,k\in K.
\end{equation}
Denote by $\mathcal H_{\chi}$ the Hilbert space of square integrable functions on $Z$ with respect to the norm  
\begin{equation}
	\small
	(f,f):=\left(\frac{i}{2}\right)^6\int_Z |f(z)|^2\,dz_{21}\wedge d\overline{z_{21}} \wedge dz_{31}\wedge d\overline{z_{31}}\wedge dz_{32}\wedge d\overline{z_{32}}\wedge dz_{41}\wedge d\overline{z_{41}} \wedge dz_{42}\wedge d\overline{z_{42}}\wedge dz_{43}\wedge d\overline{z_{43}}.%,\,\,f,h\in\mathcal H_{\chi}.
\end{equation}
Let $\left(T^{\chi}(g)f\right)(z):=\chi(k_g)\cdot\beta^{-1/2}(k_g)\cdot f(z\bar g)=$
\begin{equation}
|k_{22}|^{-2+m_2+\sqrt{-1}\rho_2}\cdot k_{22}^{-m_2}\cdot|k_{33}|^{-4+m_3+\sqrt{-1}\rho_3}\cdot k_{44}^{-m_4} \cdot|k_{44}|^{-6+m_4+\sqrt{-1}\rho_4}\cdot k_{44}^{-m_4} \cdot f(z\bar g),\,\,f\in\mathcal H_{\chi}.
\end{equation}

There are three types of degenerate principal series as follows.  Firstly, let
\begin{equation}
	K_{3,1}:=\left\{\left.k=\left(\begin{matrix}
		k_{11}&	k_{12}&	k_{13}&	k_{14}\\
		k_{21}&k_{22}&k_{23}&k_{24}\\
		k_{31}&k_{32}&k_{33}&k_{34}\\
		0&0&0&k_{44}\\
	\end{matrix}\,\right)\right|k_{ij}\in\mathbb C,\,\,k_{44}\cdot \begin{vmatrix}
		k_{11}&k_{12}&k_{13}\\
		k_{21}&k_{22}&k_{23}\\
		k_{31}&k_{32}&k_{33}\\
	\end{vmatrix}=1\right\}.
\end{equation}
Then the modular function is $\beta(k)=\frac{d\mu_l(k)}{d\mu_r(k)}=|k_{44}|^{8}$. Define
\begin{equation}
	Z_{3,1}:=\left\{\left.z=\left(\begin{matrix}
		1&0 &0&0\\
		0&1&0&0\\
		0&0&1&0\\
		z_{41}&z_{42}&z_{43}&1\\
	\end{matrix}\,\,\right)\right|z_{41},z_{42},z_{43}\in\mathbb C\,\right\}.
\end{equation}
For $m_2\in\mathbb Z$ and $\rho_2\in\mathbb R$, define a character $\chi=(m_2,\rho_2)$ of $K_{3,1}$ by
\begin{equation}
	\chi(k):=|k_{44}|^{m_2+\sqrt{-1}\rho_2}\cdot k_{44}^{-m_2},\,\,\,\,k\in K_{3,1}.
\end{equation} 
Denote by $\mathcal H_{\chi}$  the Hilbert space of square integrable functions on $Z_{3,1}$ with respect to the norm
\begin{equation}
	(f,f):=\left(\frac{i}{2}\right)^3\int_{Z_{3,1}} |f(z)|^2\,dz_{41}\wedge d\overline{z_{41}}\wedge dz_{42}\wedge d\overline{z_{42}}\wedge dz_{43}\wedge d\overline{z_{43}}.%,\,\,\,\,f,h\in\mathcal H_{\chi}.
\end{equation}
Decomposing $z\cdot g=k_g\cdot (z\bar g)$, we can define
\begin{equation}
	\left(T^{\chi}(g)f\right)(z):=\chi(k_g)\cdot\beta^{-1/2}(k_g)\cdot f(z\bar g)=|k_{44}|^{-4+m_2+\sqrt{-1}\rho_2}\cdot k_{44}^{-m_2}\cdot f(z\bar g),\,\,f\in\mathcal H_{\chi},
\end{equation}
where $k_{44}$, $z_{ij}^{\prime}$ are defined by (\ref{k44}), (\ref{z41}), respectively. 

Secondly, consider
\begin{equation}
K_{2,1,1}:=\left\{\left.k=\left(\begin{matrix}
k_{11}&	k_{12}&	k_{13}&	k_{14}\\
k_{21}&k_{22}&k_{23}&k_{24}\\
0&0&k_{33}&k_{34}\\
0&0&0&k_{44}\\
	\end{matrix}\,\right)\right|k_{ij}\in\mathbb C,\,\,k_{33}\cdot k_{44}\cdot \begin{vmatrix}
	k_{11}&k_{12}\\
	k_{21}&k_{22}\\
	\end{vmatrix}=1\right\}.
\end{equation}
The modular function is $\beta(k):=|k_{33}|^{6}\cdot|k_{44}|^{10}$. Define
\begin{equation}
	Z_{2,1,1}:=\left\{\left.z=\left(\begin{matrix}
		1&0 &0&0\\
		0&1&0&0\\
		z_{31}&z_{32}&1&0\\
		z_{41}&z_{42}&z_{43}&1\\
	\end{matrix}\,\,\right)\right|	z_{31},z_{32},z_{41},z_{42},z_{43}\in\mathbb C\,\right\}.
\end{equation}
For $m_2,m_3\in\mathbb Z$ and $\rho_2,\rho_3\in\mathbb R$, define a character $\chi=(m_2,m_3,\rho_2,\rho_3)$ of $K_{2,1,1}$ by
\begin{equation}
	\chi(k):=|k_{33}|^{m_2+\sqrt{-1}\rho_2}\cdot k_{33}^{-m_2}\cdot|k_{44}|^{m_3+\sqrt{-1}\rho_3}\cdot k_{44}^{-m_3},\,\,\,\,k\in K_{2,1,1}.
\end{equation} 
Let $\mathcal H_{\chi}$ be the Hilbert space of square integrable functions on $Z_{2,1,1}$ with respect to the norm
\begin{equation}
	(f,f):=\left(\frac{i}{2}\right)^5\int_{Z_{2,1,1}} |f(z)|^2\,dz_{31}\wedge d\overline{z_{31}}\wedge dz_{32}\wedge d\overline{z_{32}}\wedge dz_{41}\wedge d\overline{z_{41}}\wedge dz_{42}\wedge d\overline{z_{42}}\wedge dz_{43}\wedge d\overline{z_{43}}.%,\,\,\,\,f,h\in\mathcal H_{\chi}.
\end{equation}
Decomposing $z\cdot g=k_g\cdot (z\bar g)$, we can define
\begin{equation}
\begin{split}
\left(T^{\chi}(g)f\right)(z):&=\chi(k_g)\cdot\beta^{-1/2}(k_g)\cdot f(z\bar g)\\
&=|k_{33}|^{-3+m_2+\sqrt{-1}\rho_2}\cdot k_{33}^{-m_2}\cdot|k_{44}|^{-5+m_3+\sqrt{-1}\rho_3}\cdot k_{44}^{-m_3},\,\,f\in\mathcal H_{\chi}.
\end{split}
\end{equation}
Here $k_{33}$, $k_{44}$, $z_{31}^{\prime}$, $z_{32}^{\prime}$ are defined by (\ref{k33}), (\ref{k44}), (\ref{z31}), (\ref{z32}),  respectively; $z_{41}^{\prime},z_{42}^{\prime},z_{43}^{\prime}$ are defined by (\ref{z41}).

Thirdly, consider
\begin{equation}\label{K22}
	K_{2,2}:=\left\{\left.k=\left(\begin{matrix}
		k_{11}&	k_{12}&	k_{13}&	k_{14}\\
		k_{21}&k_{22}&k_{23}&k_{24}\\
		0&0&k_{33}&k_{34}\\
		0&0&k_{43}&k_{44}\\
	\end{matrix}\,\right)\right|k_{ij}\in\mathbb C,\,\,\begin{vmatrix}
	k_{11}&k_{12}\\
	k_{21}&k_{22}\\
\end{vmatrix}\cdot \begin{vmatrix}
		k_{33}&k_{34}\\
		k_{43}&k_{44}\\
	\end{vmatrix}=1\right\}.
\end{equation}
The modular function is $\beta(k):=\left|
k_{33}\cdot k_{44} -k_{34}\cdot k_{43} \right|^8$ (see \S 11 of \cite{GN2}). Define
\begin{equation}\label{Z22}
	Z_{2,2}:=\left\{\left.z=\left(\begin{matrix}
		1&0 &0&0\\
		0&1&0&0\\
		z_{31}&z_{32}&1&0\\
		z_{41}&z_{42}&0&1\\
	\end{matrix}\,\,\right)\right|	z_{31},z_{32},z_{41},z_{42}\in\mathbb C\,\right\}.
\end{equation}
For $m_2\in\mathbb Z$ and $\rho_2\in\mathbb R$, define a character $\chi=(m_2,\rho_2)$ of $K_{2,2}$ by
\begin{equation}
	\chi(k):=\left|
	k_{33}\cdot k_{44} -k_{34}\cdot k_{43} \right|^{m_2+\sqrt{-1}\rho_2}\cdot \left(k_{33}\cdot k_{44} -k_{34}\cdot k_{43}\right)^{-m_2},\,\,\,\,k\in K_{2,2}.
\end{equation} 
Let $\mathcal H_{\chi}$ be the Hilbert space of square integrable functions on $Z_{2,2}$ with respect to the norm
\begin{equation}
	(f,f):=\left(\frac{i}{2}\right)^4\int_{Z_{2,2}} |f(z)|^2\,dz_{31}\wedge d\overline{z_{31}}\wedge dz_{32}\wedge d\overline{z_{32}}\wedge dz_{41}\wedge d\overline{z_{41}}\wedge dz_{42}\wedge d\overline{z_{42}}.%,\,\,\,\,f,h\in\mathcal H_{\chi}.
\end{equation}
Decomposing $z\cdot g=k_g\cdot (z\bar g)$, we can define
\begin{equation}
	\begin{split}
		\left(T^{\chi}(g)f\right)&(z):=\chi(k_g)\cdot\beta^{-1/2}(k_g)\cdot f(z\bar g)\\
		&=\left|
		k_{33}\cdot k_{44} -k_{34}\cdot k_{43} \right|^{-4+m_2+\sqrt{-1}\rho_2}\cdot \left(k_{33}\cdot k_{44} -k_{34}\cdot k_{43}\right)^{-m_2}\cdot f(z\bar g),\,\,f\in\mathcal H_{\chi}.
	\end{split}
\end{equation}
Here $z_{31}^{\prime}$, $z_{32}^{\prime}$, $z_{41}^{\prime}$, $z_{42}^{\prime}$ are defined by  (\ref{z31})-(\ref{z41}),  and $k_{33}\cdot k_{44} -k_{34}\cdot k_{43}=$
\begin{equation}\label{k33k44}
	\footnotesize
\begin{vmatrix}
	g_{33}&g_{34}\\
	g_{43}&g_{44}\\
\end{vmatrix}+\begin{vmatrix}
	g_{32}&g_{34}\\
	g_{42}&g_{44}\\
\end{vmatrix} z_{32}+\begin{vmatrix}
	g_{31}&g_{34}\\
	g_{41}&g_{44}\\
\end{vmatrix} z_{31}+\begin{vmatrix}
	g_{32}&g_{33}\\
	g_{42}&g_{43}\\
\end{vmatrix}\begin{vmatrix}
	z_{32}&1\\
	z_{42}&z_{43}\\
\end{vmatrix}+\begin{vmatrix}
	g_{31}&g_{33}\\
	g_{41}&g_{43}\\
\end{vmatrix}\begin{vmatrix}
	z_{31}&1\\
	z_{41}&z_{43}\\
\end{vmatrix}+\begin{vmatrix}
	g_{31}&g_{32}\\
	g_{41}&g_{42}\\
\end{vmatrix}\begin{vmatrix}
	z_{31}&z_{32}\\
	z_{41}&z_{42}\\
\end{vmatrix}.
\end{equation}

\smallskip

In the following, we shall describe four types of (degenerate) complementary series. Firstly, consider the upper triangular group $K$ of $SL_4(\mathbb{C})$. For  $m_2,m_3\in\mathbb Z$, $\rho_2,\rho_3\in\mathbb R$, and $0<\sigma<1$, define a character $\chi=(m_2,m_3,\rho_2,\rho_3,\sigma)$ by 
\begin{equation}
\chi(k):=|k_{22}|^{m_2+\sqrt{-1}\rho_2}\cdot k_{22}^{-m_2}\cdot|k_{33}|^{m_3+\sqrt{-1}\rho_3+\sigma}\cdot k_{33}^{-m_3}\cdot |k_{44}|^{m_3+\sqrt{-1}\rho_3-\sigma}\cdot k_{44}^{-m_3}, \,\,k\in K.
\end{equation}
Define unipotent subgroups 
\begin{equation}
%	\footnotesize
	\dot{Z}:=\left\{\left.\dot{z}=\left(\begin{matrix}
		1&0 &0&0\\
		0&1&0&0\\
		0&0&1&0\\
		0&0&z_{1}&1\\
	\end{matrix}\,\,\right)\right|z_{1}\in\mathbb C\,\right\},\,\,\,	\widehat Z:=\left\{\left. \hat{z}= \left(\begin{matrix}
	1&0 &0&0\\
	z_{21}&1&0&0\\
	z_{31}&z_{32}&1&0\\
	z_{41}&z_{42}&0&1\\
	\end{matrix}\,\,\right)\right|z_{ij}\in\mathbb C\,\right\},
\end{equation}
and write $z=(\hat{z},z_1)$.  Let $d\mu:=dz_{21}\wedge d\overline{z_{21}}\wedge dz_{31}\wedge d\overline{z_{31}}\wedge dz_{41}\wedge d\overline{z_{41}}\wedge dz_{32}\wedge d\overline{z_{32}}\wedge dz_{42}\wedge d\overline{z_{42}}$. Define $\mathcal H_{\chi}$ to be the Hilbert space of square integrable functions on $Z$ with respect to the norm
\begin{equation}
(f,f):=\left(\frac{i}{2}\right)^7\int_{\mathbb C\times Z} \frac{f(\hat{z},z_1)\cdot\overline{f(\hat{z},w_1)}}{\left|z_1-w_1\right|^{2-2\sigma}}\,dw_{1}\wedge d\overline{w_{1}}\wedge dz_{1}\wedge d\overline{z_{1}}\wedge d\mu. %,\,\,f,h\in\mathcal H_{\sigma}.
\end{equation}
Set $\left(T^{\chi}(g)f\right)(z):=\chi(k_g)\cdot\beta^{-1/2}(k_g)\cdot f(z\bar g)=$
\begin{equation}
|k_{22}|^{-2+m_2+\sqrt{-1}\rho_2}\cdot k_{22}^{-m_2}\cdot|k_{33}|^{-4+m_3+\sqrt{-1}\rho_3+\sigma}\cdot k_{33}^{-m_3}\cdot |k_{44}|^{-6+m_3+\sqrt{-1}\rho_3-\sigma}\cdot k_{44}^{-m_3}\cdot f(z\bar g).
\end{equation}

Secondly, we still consider $K$. For  $m_2\in\mathbb Z$, $\rho_2\in\mathbb R$, and $0<\sigma_1,\sigma_2<1$, define a character $\chi=(m_2,\rho_2,\sigma)$ by 
\begin{equation}
\chi(k):=|k_{22}|^{-2\sigma_1}\cdot|k_{33}|^{m_2+\sqrt{-1}\rho_2+\sigma_2-\sigma_1}\cdot k_{33}^{-m_2}\cdot |k_{44}|^{m_2+\sqrt{-1}\rho_2-\sigma_2-\sigma_1}\cdot k_{44}^{-m_2}, \,\,k\in K.
\end{equation}
Define unipotent subgroups 
\begin{equation}
	%	\footnotesize
	\dot{Z}:=\left\{\left.\dot{z}=\left(\begin{matrix}
		1&0 &0&0\\
		z_2&1&0&0\\
		0&0&1&0\\
		0&0&z_{1}&1\\
	\end{matrix}\,\,\right)\right|z_{1}\in\mathbb C\,\right\},\,\,\,	\widehat Z:=\left\{\left. \hat{z}= \left(\begin{matrix}
		1&0 &0&0\\
		0&1&0&0\\
		z_{31}&z_{32}&1&0\\
		z_{41}&z_{42}&0&1\\
	\end{matrix}\,\,\right)\right|z_{ij}\in\mathbb C\,\right\},
\end{equation}
and write $z=(\hat{z},z_1,z_2)$.  Let $d\mu:=dz_{31}\wedge d\overline{z_{31}}\wedge dz_{41}\wedge d\overline{z_{41}}\wedge dz_{32}\wedge d\overline{z_{32}}\wedge dz_{42}\wedge d\overline{z_{42}}$. Define $\mathcal H_{\chi}$ to be the Hilbert space of square integrable functions on $Z$ with respect to the norm
\begin{equation}
	\small
	(f,f):=\left(\frac{i}{2}\right)^6\int_{\mathbb C\times\mathbb C\times Z} \frac{f(\hat{z},z_1,z_2)\cdot\overline{f(\hat{z},w_1,w_2)}}{\left|z_1-w_1\right|^{2-2\sigma_1}\left|z_2-w_2\right|^{2-2\sigma_2}}\,dw_{1}\wedge d\overline{w_{1}}\wedge dz_{1}\wedge d\overline{z_{1}}\wedge dw_{2}\wedge d\overline{w_{2}}\wedge dz_{2}\wedge d\overline{z_{2}}\wedge d\mu. %,\,\,f,h\in\mathcal H_{\sigma}.
\end{equation}
Set $\left(T^{\chi}(g)f\right)(z):=\chi(k_g)\cdot\beta^{-1/2}(k_g)\cdot f(z\bar g)=$
\begin{equation}
|k_{22}|^{-2\sigma_1}\cdot|k_{33}|^{-4+m_2+\sqrt{-1}\rho_2+\sigma_2-\sigma_1}\cdot k_{33}^{-m_2}\cdot |k_{44}|^{-6+m_2+\sqrt{-1}\rho_2-\sigma_2-\sigma_1}\cdot k_{44}^{-m_2}\cdot f(z\bar g).
\end{equation}

Thirdly, consider the group $K_{2,2}$ defined in (\ref{K22}). For  $m_2\in\mathbb Z$, $\rho_2\in\mathbb R$, and $0<\sigma<1$, define a character $\chi=(m_2,\rho_2,\sigma)$ by 
\begin{equation}
\chi(k):=|k_{33}|^{m_2+\sqrt{-1}\rho_2+\sigma}\cdot k_{33}^{-m_2}\cdot |k_{44}|^{m_2+\sqrt{-1}\rho_2-\sigma}\cdot k_{44}^{-m_2}, \,\,k\in K_{2,2}.
\end{equation}
Define unipotent subgroups 
\begin{equation}
	%	\footnotesize
	\dot{Z}:=\left\{\left.\dot{z}=\left(\begin{matrix}
		1&0 &0&0\\
		0&1&0&0\\
		0&0&1&0\\
		0&0&z_{1}&1\\
	\end{matrix}\,\,\right)\right|z_{1}\in\mathbb C\,\right\},\,\,\,	\widehat Z:=\left\{\left. \hat{z}= \left(\begin{matrix}
		1&0 &0&0\\
		0&1&0&0\\
		z_{31}&z_{32}&1&0\\
		z_{41}&z_{42}&0&1\\
	\end{matrix}\,\,\right)\right|z_{ij}\in\mathbb C\,\right\},
\end{equation}
and write $z=(\hat{z},z_1)$.  Let $d\mu:=dz_{31}\wedge d\overline{z_{31}}\wedge dz_{41}\wedge d\overline{z_{41}}\wedge dz_{32}\wedge d\overline{z_{32}}\wedge dz_{42}\wedge d\overline{z_{42}}$. Define $\mathcal H_{\chi}$ to be the Hilbert space of square integrable functions on $Z$ with respect to the norm
\begin{equation}
	(f,f):=\left(\frac{i}{2}\right)^6\int_{\mathbb C\times Z} \frac{f(\hat{z},z_1)\cdot\overline{f(\hat{z},w_1)}}{\left|z_1-w_1\right|^{2-2\sigma}}\,dw_{1}\wedge d\overline{w_{1}}\wedge dz_{1}\wedge d\overline{z_{1}}\wedge d\mu. %,\,\,f,h\in\mathcal H_{\sigma}.
\end{equation}
Set $\left(T^{\chi}(g)f\right)(z):=\chi(k_g)\cdot\beta^{-1/2}(k_g)\cdot f(z\bar g)=$
\begin{equation}
|k_{33}|^{-3+m_2+\sqrt{-1}\rho_2+\sigma}\cdot k_{33}^{-m_2}\cdot |k_{44}|^{-5+m_2+\sqrt{-1}\rho_2-\sigma}\cdot k_{44}^{-m_2}\cdot f(z\bar g).
\end{equation}

At last, we introduce Stein's complementary series. Consider the subgroups $K_{2,2}$, $Z_{22}$ defined in (\ref{K22}), (\ref{Z22}), respectively. %Each $g=\left(\begin{matrix}A&B\\C&D\end{matrix}\right)\in SL_4(\mathbb C)$, where $A,B,C,D$ are $2\times 2$ block matrices,  acts on $z=\left(\begin{matrix}I_{2\times2}&0\\\widetilde z&I_{2\times2}\end{matrix}\right)\in SL_4(\mathbb C)$ by $\widetilde z\dashmapsto(A+\widetilde zC)^{-1}(B+\widetilde zD)$.
For   $0<\sigma<1$, define a character $\chi=(\sigma)$ by 
\begin{equation}
\chi(k):=\left|
k_{33}\cdot k_{44} -k_{34}\cdot k_{43} \right|^{-2\sigma}, \,\,k\in K_{2,2}.
\end{equation}
Let $\mathcal H_{\chi}$ be the Hilbert space of square integrable functions on $Z_{2,2}$ with respect to the norm
\begin{equation}
	(f,f):=\left(\frac{i}{2}\right)^8\int_{Z_{2,2}\times Z_{2,2}} \frac{f(z)\overline{f(w)}}{\left|\det\left(\begin{matrix}z_{31}-w_{31}&z_{32}-w_{32}\\
	z_{41}-w_{41}&z_{42}-w_{42}\\		\end{matrix}\right)\right|^{4-2\sigma}}\,d\mu\wedge d\nu, %,\,\,\,\,f,h\in\mathcal H_{\chi}.
\end{equation}
where
\begin{equation}
\begin{split}
&d\mu:=dz_{31}\wedge d\overline{z_{31}}\wedge dz_{32}\wedge d\overline{z_{32}}\wedge dz_{41}\wedge d\overline{z_{41}}\wedge dz_{42}\wedge d\overline{z_{42}},\\
&d\nu:=dw_{31}\wedge d\overline{w_{31}}\wedge dw_{32}\wedge d\overline{w_{32}}\wedge dw_{41}\wedge d\overline{w_{41}}\wedge dw_{42}\wedge d\overline{w_{42}}.
\end{split}
\end{equation}
Decomposing $z\cdot g=k_g\cdot (z\bar g)$, we can define
\begin{equation}	\left(T^{\chi}(g)f\right)(z):=\chi(k_g)\cdot\beta^{-1/2}(k_g)\cdot f(z\bar g)=\left|k_{33}\cdot k_{44} -k_{34}\cdot k_{43}\right|^{-4-2\sigma}\cdot f(z\bar g),\,\,f\in\mathcal H_{\chi}.
\end{equation}
Here $z_{31}^{\prime}$, $z_{32}^{\prime}$, $z_{41}^{\prime}$, $z_{42}^{\prime}$ are defined by  (\ref{z31})-(\ref{z41}); $k_{33}\cdot k_{44} -k_{34}\cdot k_{43}$ is given by (\ref{k33k44}).

\end{document}